\documentclass[12pt]{amsart}
\usepackage{enumerate,color}
\usepackage{graphicx}
\usepackage[text={6in,9in},centering]{geometry}
 
 \usepackage{pgfplots}
\pgfplotsset{compat=1.18}

\usepackage{tikz}
\usetikzlibrary{positioning, decorations.pathreplacing, calc}
\usepackage{amsmath, amssymb}

\theoremstyle{plain}
\newtheorem{thm}{Theorem}[section]
\newtheorem{cor}[thm]{Corollary}
\newtheorem{lem}[thm]{Lemma}
\newtheorem{prop}[thm]{Proposition}
\theoremstyle{definition}

\theoremstyle{remark}

\theoremstyle{plain}
\newtheorem{mainthm}{Theorem}

\newcommand{\bu}{\mathbf{u}}

\title[]{Hausdorff and Packing Measures for $\beta$-L{\"u}roth Infinite Iterated Function Systems}

\date{} 

\author[J. J. Miao]{Jun Jie Miao}
\thanks{Corresponding author. E-mail: jjmiao@math.ecnu.edu.cn}
\address[J. J. Miao]{School of Mathematical Sciences,  Key Laboratory of MEA(Ministry of Education) \& Shanghai Key Laboratory of PMMP,  East China Normal University, Shanghai 200241, P.~R. China}
\email{jjmiao@math.ecnu.edu.cn}

\author[Tianrui Wang]{Tianrui Wang}
\address[Tianrui Wang]{School of Mathematical Sciences,  Key Laboratory of MEA(Ministry of Education) \& Shanghai Key Laboratory of PMMP,  East China Normal University, Shanghai 200241, P.~R. China}
\email{51265500036@stu.ecnu.edu.cn}

\begin{document}
\begin{abstract}
We study the Hausdorff and packing measures of attractors generated by infinite $\beta$-L\"uroth iterated function systems. Let $I\subset\mathbb N$ be infinite and let $\dim_H E_I=h$. For a regular system, we first give conditions for the positivity of $\mathcal H^h(E_I)$ and the finiteness of $\mathcal P^h(E_I)$ in terms of first-level L\"uroth intervals. We then introduce tail-regular partitions and the scale $R_n=\frac{B_n^\infty}{B_n}.$ Under tail-regularity, the Hausdorff and packing conditions reduce to counting conditions for $I$ on intervals of length comparable to $R_n$, together with a tail condition in the packing case. Finally, we apply these results to regularly varying partitions and stretched-exponential partitions.
\end{abstract}

\maketitle

\section{Introduction}
\subsection{Infinite conformal iterated function systems.}\label{subsec_ICIFS}
A classical iterated function system (IFS) is a finite family of contractions on a complete metric space. By Hutchinson's theorem, there is a unique non-empty compact attractor $E$ satisfying
\[
E=\bigcup_{i=1}^N\varphi_i(E);
\]
see \cite{B93,E90,F14,H81}. If the contractions are conformal, the attractor is called a self-conformal set. Under the open set condition (OSC), its Hausdorff dimension is the unique zero $h$ of the pressure, $P(h)=0$; see \cite{B79,MU03,P97,PU10,R82}. If the contractions are similarities   with contraction ratios $r_i$, the attractor is called a self-similar set, which is a special case of self-conformal sets, and the dimension formula $\dim_H E=h$ is simplified into 
\[
P(t)=\log\sum_{i=1}^N r_i^t,
\qquad
\sum_{i=1}^N r_i^h=1,
\]
   under the OSC \cite{E90,F14,H81,M46}.

The OSC is also important for the study of the   measures supported on such attractors. For   self-similar sets, Bandt--Graf and Schief related the OSC to positivity of the Hausdorff measure in the similarity dimension \cite{BG92,S94}, and the corresponding self-conformal result was obtained by Peres, Rams, Simon and Solomyak \cite{PRSS01}. 
For self-conformal sets, the OSC is equivalent to the positivity of the $h$-dimensional Hausdorff measure \cite{PRSS01}. If, in addition, the self-conformal sets satisfy the cone condition, then the $h$-dimensional packing measures are also positive and finite; see \cite{MU96,TT85}.

Countably infinite conformal iterated function systems were systematically studied by Mauldin and Urbański \cite{MU96,MU99}; see also \cite{MU03}. They extended the thermodynamic formalism of the classical conformal iterated function systems to countable infinite settings and related pressures to dimension, conformal measures, and the geometry of the limit set. We first recall some basic  definitions.

Let $I$ be a countable index set and $X\subset\mathbb{R}^d$ be a compact subset with $\operatorname{int} X\neq\emptyset$. We call
$\Phi_I=\big\{\varphi_i:X\to X\big\}_{i\in I}$ \emph{an infinite conformal iterated function
system (ICIFS) } if it satisfies the following:
\begin{itemize}
\item[(i)] \emph{Uniform contraction}: there exists $0<c<1$ such that for
all $i\in I$,
\begin{equation*}
|\varphi_i(x)-\varphi_i(y)|\le c|x-y|
\qquad\text{for all }x,y\in X.
\end{equation*}
\item[(ii)]\emph{Open set  condition (OSC)}: For all distinct $i,j\in I$,
\[
\varphi_i(\operatorname{int} X)\cap \varphi_j(\operatorname{int} X)=\emptyset.
\]
\item[(iii)] \emph{Conformality}: there exists an open connected set $V$
with $X\subset V$ such that each $\varphi_i$ extends to a $C^1$
conformal diffeomorphism from $V$ into $V$.
\item[(iv)] \emph{Bounded distortion}: there exists $C\ge1$ such that,
for every word $\bu=u_1\ldots u_k\in I^k$, the composition
$\varphi_\bu=\varphi_{u_1}\circ\cdots\circ\varphi_{u_k}$ satisfies
\[
\|\varphi_\bu'(x)\|\le C\|\varphi_\bu'(y)\|  \qquad \qquad  \textit{
for  $x,y\in V$.}
\]
\item[(v)] \emph{Cone condition}: the initial compact set $X$ satisfies
the cone condition, namely, there exist $\alpha,l>0$ such that for every
$x\in\partial X$, there is an open cone
\[
\operatorname{Con}(x,u_x,\alpha,l)\subset \operatorname{int}(X)
\]
with vertex $x$, direction vector $u_x$, central angle $\alpha$, and
altitude $l$.
\end{itemize}
For $\omega=(\omega_1,\omega_2,\ldots)\in I^{\mathbb N}$, uniform
contractivity implies that
\[
\bigcap_{n=1}^{\infty}\varphi_{\omega|_n}(X)
\]
consists of a single point, denoted by $\pi(\omega)$. The limit set is $E=E(\Phi_I)=\pi(I^{\mathbb N})$.

Following Mauldin and Urbański \cite{MU96,MU99}, for $n\ge1$ set
\[
\psi_n(t)=\sum_{\omega\in I^n}\|\varphi_\omega'\|^t,
\]
where $\|\varphi_\omega'\|=\sup_{x\in X}\|\varphi_\omega'(x)\|$, and define
\[
\theta_{\Phi_I}=\inf\big\{t:\psi_1(t)<\infty\big\}.
\]
The topological pressure is
\[
P(t)=\lim_{n\to\infty}\frac1n\log\psi_n(t),
\]
and the critical exponent is
\[
h_{\Phi_I}=\inf\big\{t:P(t)<0\big\}.
\]
Mauldin and Urbański proved that $P(t)$ is non-increasing on $[0,\infty)$ and, on its finiteness domain $\{t\ge0:P(t)<\infty\}$, is strictly decreasing, continuous, and convex; moreover, $P(d)\le0$; see \cite{MU96,MU99} for details.

Given a subset $J\subset I$, we write
\[
E_J = \pi(J^{\mathbb N})
\]
for the limit set of the subsystem $\{\varphi_i\}_{i\in J}$. We call $\{\varphi_i\}_{i\in J}$ \emph{cofinite} if $I\setminus J$ is finite. If $J$ is finite,
$E_J$ is the attractor of the finite conformal iterated function system $\{\varphi_i\}_{i\in J}$. In \cite{MU96,MU99}, they proved the following conclusions.

\begin{mainthm}\label{dimH}
The Hausdorff dimension of the limit set satisfies
\[
\dim_H E=\sup\big\{\dim_H E_J:J\subset I\text{ is finite}\big\}=h_{\Phi_I}.
\]
If $P(t)=0$, then $t=\dim_H E$.
\end{mainthm}

The system $\Phi_I$ is called \emph{regular} if there exists $t$ such that
$P(t)=0$. Regularity is
equivalent to the existence of a $t$-conformal measure. A Borel probability
measure $m$ on $X$ is \emph{$t$-conformal} if $m(E)=1$, and for every Borel set
$A\subset X$ and every $i\in I$,
\begin{equation}\label{conformal measure}
m(\varphi_i(A))=\int_A |\varphi_i'|^t\,dm,
\end{equation}
with $m(\varphi_i(X)\cap\varphi_j(X))=0$
for every $i\ne j$; see \cite{MU96,MU99}.
The system $\Phi_I$ is called \emph{absolutely regular} if every non-empty subsystem is regular \cite{MU99}.
\begin{mainthm}\label{mu-absolute-regularity}
 $\Phi_I$ is absolutely regular if and only if $\theta_{\Phi_I}=0$.
\end{mainthm}

Note that for self-conformal sets,  the pressures always have a zero, however, for countably infinite systems, the regularity is not automatically satisfied. Moreover, even for regular infinite systems, the critical Hausdorff and packing measures of the limit set need not have the positivity and finiteness properties that hold in the finite case. Mauldin and Urbański developed characterizations of these measures in terms of conformal measures; related results for infinitely generated self-similar and continued-fraction-type systems were obtained in \cite{MU96,MU99,M96}.  

\subsection{The $\beta$-L{\"u}roth systems}
Continued fractions are a classical tool in Diophantine approximation and metric number theory \cite{IK02,K64,RS92}. Restricting the allowed partial quotients produces Cantor-type sets, and the study of their  dimensions, together with the associated pressure equations and measure properties, has a long history; see \cite{GM21,G41,H06,MU99}. L\"uroth expansions provide a closely related countable symbolic coding of the interval, but with affine rather than M\"obius inverse branches. Following L\"uroth's original work \cite{L83}, generalized L\"uroth expansions and $\alpha$-L\"uroth expansions have been studied from ergodic, metric and thermodynamic viewpoints \cite{BBDK96,KMS12}, as well as from the point of view of dimension theory \cite{ABGR22,BI09,TZ21}. Since the inverse branches of such expansions form  infinite self-similar systems, this makes L\"uroth-type expansions a natural setting for our study.

The questions studied in this paper are motivated by Mauldin and Urba\'nski \cite{MU99}. For restricted continued fractions, they related the critical Hausdorff and packing measures to arithmetic properties of the set of allowed digits. Fishman, Simmons and Urba\'nski later reformulated the continued-fraction condition through local counting and tail conditions \cite{FSU14}. We study the corresponding measure problem for general $\beta$-L\"uroth partitions. We begin by recalling the basic definitions.

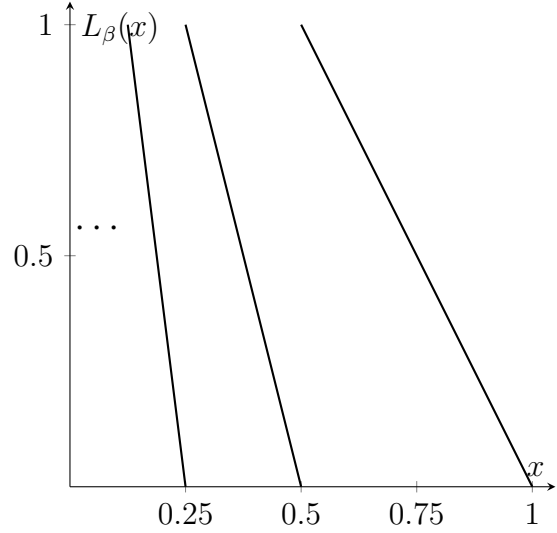
\begin{figure}[htbp]
\centering
\begin{tikzpicture}
\begin{axis}[
    width=8cm,
    height=8cm,
    xlabel={$x$},
    ylabel={$L_\beta(x)$},
    xmin=0, xmax=1.05,
    ymin=0, ymax=1.05,
    xtick={0, 0.25, 0.5, 0.75, 1},
    ytick={0, 0.5, 1},
    axis lines=middle,
    samples=2,
    clip=false,
    every axis plot/.append style={thick, black},
]

\addplot[domain=0.5:1, black] { (1 - x)/0.5 };
\addplot[domain=0.25:0.5, black] { (0.5 - x)/0.25 };
\addplot[domain=0.125:0.25, black] { (0.25 - x)/0.125 };
\node[scale=1.2] at (axis cs:0.065,0.56) {$\cdots$};

\end{axis}
\end{tikzpicture}
\caption{Example of the $\beta$-Lüroth map for the    partition $B_n = 2^{-n}$. }
\label{fig:beta_luroth} \label{fig_Lurothex}
\end{figure}

 Let $\beta=\big\{\mathbf{B}_n\big\}_{n=1}^\infty$ be a countably infinite partition
of $(0,1]$ into non-empty, right-closed and left-open intervals. Write
\begin{equation*}
B_n=\lambda(\mathbf{B}_n),
\qquad
B_k^l=\sum_{n=k}^l B_n,
\qquad
B_k^\infty=\sum_{n=k}^\infty B_n.
\end{equation*}
We adopt the convention that \(B_k^l=0\) whenever \(k>l\). We always assume that the intervals are ordered from right to left,
starting from $\mathbf{B}_1$, and that $B_{n+1}< B_n$ for all
$n\in\mathbb N$.  We define the $\beta$-L{\"u}roth map
$L_\beta:[0,1]\to[0,1]$ by
\[
L_\beta(x)=
\begin{cases}
\frac{B_n^\infty-x}{B_n}, & x\in\mathbf{B}_n,\ n\in\mathbb N,\\[1.2ex]
0, & x=0.
\end{cases}
\]
See Figure \ref{fig_Lurothex} for an example. The map $L_\beta$ generates a L{\"u}roth-type expansion of points in $[0,1]$. For $x\in(0,1]$, the digits are defined successively for as long as
$L_\beta^{k-1}(x)\neq0$. If the orbit never reaches $0$, then
\[
x
=
\sum_{n=1}^{\infty}
(-1)^{n-1}
\Bigl(\prod_{i<n}B_{\ell_i}\Bigr)B_{\ell_n}^{\infty}
=
B_{\ell_1}^{\infty}
-B_{\ell_1}B_{\ell_2}^{\infty}
+B_{\ell_1}B_{\ell_2}B_{\ell_3}^{\infty}
-\cdots,
\]
and we write $x=[\ell_1,\ell_2,\ell_3,\ldots]_\beta$.
If the expansion terminates after $k$ digits, we write
$x=[\ell_1,\ldots,\ell_k]_\beta$. Points with finite expansion are called
\emph{$\beta$-rational}, and the others \emph{$\beta$-irrational}.

Let $X=[0,1]$. For each $n\in\mathbb N$, we write
\begin{equation*}
\varphi_n(x)=B_n^\infty-B_nx, \qquad x\in X=[0,1].
\end{equation*}
For $0\leq x<1$, $\varphi_n$ is the inverse of the $n$-th branch of
$L_\beta$; at $x=1$ it is the affine extension to the closure, with
$\varphi_n(1)=B_{n+1}^\infty$, the left endpoint of $\mathbf B_n$.
Given an infinite set $I\subset\mathbb{N}$,
$\Phi_I=\big\{\varphi_i\big\}_{i\in I}$ is an infinite self-similar
iterated function system on $X=[0,1]$, and we call it \emph{a
$\beta$-L{\"u}roth system with index set $I$}. It is clear that $\Phi_I$
satisfies the OSC, and the cone condition is automatically satisfied.
We denote its attractor by $E_I$. When the index set is fixed, we abbreviate
$\Phi_I$ and $E_I$ by $\Phi$ and $E$. It is clear that 
$E_I=\pi(I^{\mathbb N})$. 

For a $\beta$-L{\"u}roth system with index set $I$, we have
\begin{equation}\label{def_phiIt}
\psi_1(t)=\sum_{n\in I}(B_n)^t,
\qquad
P(t)=\log\psi_1,
\qquad
\theta_I=\inf\Bigl\{t:\sum_{n\in I}(B_n)^t<\infty\Bigr\}.
\end{equation}
If $\Phi_I$ is regular and $\dim_H E_I=h$, then
\begin{equation}\label{regular-pressure-eq}
P(h)=0,
\qquad
\sum_{n\in I}(B_n)^h=1.
\end{equation}

For a set $A\subset\mathbb R$
and $r>0$, we denote the open $r$-neighbourhood of $A$ by
\[
B(A,r)=\big\{x\in\mathbb R:\operatorname{dist}(x,A)<r\big\}.
\]
We   use  $A\lesssim B$ to mean that there exists a constant $C>0$, independent of the variables, such that $A\leq CB$; similarly, $A\gtrsim B$ means $B\lesssim A$, and $A\asymp B$ means that both $A\lesssim B$ and $A\gtrsim B$ hold.

\subsection{Main conclusions}
The main purpose of this paper is to describe absolute regularity and the critical Hausdorff and packing measures through the geometry of the L\"uroth partition and the distribution of the admissible indices. Throughout the paper, write
\[
I=\{i_1<i_2<i_3<\cdots\}.
\]
Set
\[
\mathcal S(I)=\{(B_i)^{-2}:i\in I\}.
\]
We say that $\mathcal S(I)$ has \emph{strong arithmetic density zero} if, for every $t>0$,
\[
\lim_{x\to\infty}
\frac{\#\bigl(\mathcal S(I)\cap[1,x]\bigr)}{x^t}=0.
\]

The following theorem characterizes absolute regularity of the $\beta$-L{\"u}roth system.
\begin{thm}\label{absolute-regularity}
Let $\Phi_I$ be a $\beta$-L{\"u}roth system. Then the following   are equivalent:
\begin{enumerate}
\item[(1)] $\mathcal S(I)$ has strong arithmetic density zero.
\item[(2)] For each $t>0$, $\displaystyle \sum_{n\in I}(B_n)^t<\infty.$

\item[(3)] For each $t>0$,$\displaystyle \lim_{k\to\infty}k(B_{i_k})^t=0.$
\item[(4)] $\Phi_I$ is absolutely regular.
\end{enumerate}
\end{thm}

We say the  $\beta$-L{\"u}roth system $\Phi_I$ satisfies the geometric ratio condition (GRC) if
\begin{equation}\label{hausdorff-hyp-gap}
\sup_{j\geq1}
\frac{B_{i_j}}{B_{i_j+1}^{i_{j+1}}}<\infty.
\end{equation}
Next, we give a necessary and sufficient condition for the  $h$-dimensional Hausdorff measure to be zero.

\begin{thm}\label{hm}
Given a regular $\beta$-L{\"u}roth system $\Phi_I$ satisfying GRC, let $E_I$ be the attractor of $\Phi_I$ with $\dim_H E_I=h$. Then $\mathcal H^h(E_I)=0$
if and only if
\[
\sup_{k<l}
\frac{1}
{(B_k^l)^h} \sum_{n=k}^{l}\mathbf{1}_I(n)(B_n)^h
=\infty.
\]
\end{thm}

It is natural to study the analogous question for packing measure. Since the packing case   is governed by local mass estimates at all scales, the criterion is necessarily scale-by-scale rather than a single global block estimate. 
\begin{thm}\label{packing}
Given a regular $\beta$-L{\"u}roth system $\Phi_I$ satisfying GRC, let $E_I$ be the attractor of $\Phi_I$ with $\dim_H E_I=h$.  Then $\mathcal P^h(E_I)<\infty$ if and only if there exist constants $\kappa>0$, $\xi>0$, and a finite set $J\subset I$ such that, for every $i\in I\setminus J$ and every $r$ satisfying $B_i<r\leq\xi,$
there exist indices $k\leq i\leq l$ with $k\in\mathbb N$ and $l\in\mathbb N\cup\{\infty\}$ such that
\begin{equation}\label{packing-main-geometry}
B_k^{i-1}<r,
\qquad
B_{i+1}^{l}<r,
\qquad
B_k^l<2r,
\end{equation}
and
\begin{equation}\label{packing-main-mass}
\sum_{n=k}^{l}\mathbf{1}_I(n)(B_n)^h\geq \kappa r^h.
\end{equation}
\end{thm}

The following corollary is an immediate consequence of Theorem~\ref{packing}, and it plays a key role in the proof of Theorem~\ref{tail-packing}. 
\begin{cor}\label{cor_PT}
Given a regular $\beta$-L{\"u}roth system $\Phi_I$ satisfying GRC, let $E_I$ be the attractor of $\Phi_I$ with $\dim_H E_I=h$. If
$\mathcal P^h(E_I)<\infty$, then
\[
\inf_{k\geq1}
\frac{1}{(B_k^\infty)^h}
\sum_{n=k}^{\infty}\mathbf{1}_I(n)(B_n)^h>0.
\]
\end{cor}

In general, the conditions in Theorem~\ref{hm} and Theorem~\ref{packing} are difficult to verify directly. To make them more applicable, we introduce the following tail-regularity condition on the partitions.  Given    $\beta=\big\{\mathbf{B}_n\big\}_{n=1}^\infty$,  we write
\[
R_n=\frac{B_n^\infty}{B_n}.
\]
We call $\beta$ \emph{tail-regular} if $\lim_{n\to\infty} R_n=\infty$ and there exist constants   $\eta>0$ and $n_0\in\mathbb N$ such that for all $n\ge n_0$, we have   
\[
\eta R_n\le\frac n2  \quad \textit{  and } \quad B_m\asymp  B_n
\qquad\text{ for all $m\in\mathbb{N}$ with }|m-n|\le\eta R_n, 
\]
where the constant in $\asymp$ is independent of $m$ and $n$.                                        

\begin{thm}\label{tail-hausdorff}
Suppose  $\beta$ is  tail-regular. Let $\Phi_I$ be a regular $\beta$-L{\"u}roth system   with $\dim_H E_I=h$. Then there exists $\delta>0$, depending only on $\beta$, such that $\mathcal H^h(E_I)>0$ if and only if there exists $C>0$ such that, for every sufficiently large $n$ and every integer $1\le m\le\delta R_n$,
\begin{equation}\label{tail-H-counting-intro}
\#\bigl(I\cap[n,n+m-1]\bigr)\le Cm^h.
\end{equation}
\end{thm}

\begin{thm}\label{tail-packing}
Given a tail-regular $\beta$, let $\Phi_I$ be a regular $\beta$-L{\"u}roth system   with $\dim_H E_I=h$. Then there exists $\delta>0$, depending only on $\beta$, such that $\mathcal P^h(E_I)<\infty$ if and only if the following  hold:
\begin{enumerate}
\item[(i)] there exists $c>0$ such that, for every sufficiently large $i\in I$ and every integer $1\le m\le\delta R_i$,
\begin{equation*}
\#\bigl(I\cap[i-m,i+m]\bigr)\ge cm^h;
\end{equation*}
\item[(ii)] there exists $c'>0$ such that, for every sufficiently large $n$,
\begin{equation}\label{tail-P-tail-intro}
\sum_{\substack{j\in I\\j\ge n}}(B_j)^h\ge c'(B_n^\infty)^h.
\end{equation}
\end{enumerate}
\end{thm}

The rest of the paper is organised as follows. In Section~2, we give the proof of  Theorem~\ref{absolute-regularity} which   characterizes  the absolute regularity.  We study the properties of the $h$-dimensional Hausdorff measures and give the proof of Theorem \ref{hm} in Section \ref{sec_HMD}. In Section~4, we study the finiteness of the $h$-dimensional packing measure and prove Theorem \ref{packing} and corollary \ref{cor_PT}. In the final Section, we study the properties of  tail-regularity for the  partitions of regular $\beta$-L{\"u}roth systems, and we provide the proofs for Theorems~\ref{tail-hausdorff} and \ref{tail-packing} and apply these results to regularly varying and stretched-exponential partitions.

\section{Absolute regularity of  $\beta$-L{\"u}roth systems}
In this section,  we prove Theorem~\ref{absolute-regularity}. Throughout, 
$\Phi_I=\{\varphi_i\}_{i\in I}$ is the $\beta$-L\"uroth system associated 
with an infinite index set $I=\{i_1<i_2<\cdots\}\subset\mathbb N$, and 
$\beta=\{\mathbf B_n\}_{n=1}^\infty$ is the underlying partition of $(0,1]$. 
Recall that $B_n=\lambda(\mathbf B_n),$ $B_k^l=\sum_{n=k}^l B_n$ and $B_k^\infty=\sum_{n=k}^\infty B_n,$
and that $\mathcal S(I)=\{(B_i)^{-2}:i\in I\}.$ 
Theorem~\ref{absolute-regularity} states that the absolute regularity of 
$\Phi_I$ is equivalent to each of the following: $\mathcal S(I)$ has 
strong arithmetic density zero; $\sum_{i\in I}(B_i)^t<\infty$ for every 
$t>0$; and $k(B_{i_k})^t\to0$ for every $t>0$.

\begin{proof}[Proof of Theorem~\ref{absolute-regularity}]
First, we prove that (1)$\Longleftrightarrow$(2). Suppose that $\mathcal S(I)$ has strong arithmetic density zero, that is 
\[
\lim_{x\to\infty}
\frac{\#\bigl(\mathcal S(I)\cap[1,x]\bigr)}{x^\varepsilon}=0
\qquad\text{for every }\varepsilon>0.
\]
Fix $t>0$ and set $\varepsilon=t/2$. There exists $j_0$ such that for every $j\ge j_0$, $\#\bigl(\mathcal S(I)\cap[1,2^{j+1}]\bigr)\le 2^{(j+1)\varepsilon}.$ There are only finitely many $n\in I$ with $B_n^{-2}<2^{j_0}$. Hence
\[
\begin{aligned}
\sum_{\substack{n\in I\\ (B_n)^{-2}\ge2^{j_0}}}(B_n)^t
&=\sum_{j\ge j_0}
\sum_{\substack{n\in I\\2^j\le (B_n)^{-2}<2^{j+1}}}(B_n)^t\\
&\le\sum_{j\ge j_0}2^{-jt}\#\bigl(\mathcal S(I)\cap[1,2^{j+1}]\bigr)\\
&\le 2^{\varepsilon}\sum_{j\ge j_0}2^{-j(t-\varepsilon)}<\infty.
\end{aligned}
\]
Thus $\sum_{n\in I}(B_n)^t<\infty$, i.e.,  (2) holds.

Conversely, suppose that (2) holds. For every $t>0$,
\[
\begin{aligned}
\frac{\#\bigl(\mathcal S(I)\cap[1,x]\bigr)}{x^t}
&=\sum_{\substack{(B_n)^{-2}\le x}}\frac{\mathbf{1}_I(n)}{x^t}\le x^{-t/2}\sum_{\substack{(B_n)^{-2}\le x}}\mathbf{1}_I(n)(B_n)^{t}\le x^{-t/2}\sum_{n\in I}(B_n)^{t}.
\end{aligned}
\]
Since $\sum_{n\in I}(B_n)^{t}<\infty$, it follows that
\[
\lim_{x\to\infty}
\frac{\#\bigl(\mathcal S(I)\cap[1,x]\bigr)}{x^t}=0.
\]
Hence (1) holds.

Next, we show  that (1)$\Longleftrightarrow$(3).
Since $\#\bigl(\mathcal S(I)\cap[1,B_{i_k}^{-2}]\bigr)=k,$  (1) $\Longrightarrow$ (3) is straightforward. Conversely, suppose that   $\lim_{k\to\infty}k(B_{i_k})^t=0$ for all $t>0$. Fix $t>0$. Since   $\#\bigl(\mathcal S(I)\cap[1,x]\bigr)=k$ for $B_{i_k}^{-2}\le x<B_{i_{k+1}}^{-2}$, we have 
\[
\lim_{x\to\infty} \frac{\#\bigl(\mathcal S(I)\cap[1,x]\bigr)}{x^t}\le \lim_{k\to\infty}kB_{i_k}^{2t} =0.
\]
Thus (3) $\Longrightarrow$ (1).

Finally, we only need to show (2)$\Longleftrightarrow$(4). Recall that
\[
\theta_I=\inf\Bigl\{t:\sum_{n\in I}(B_n)^t<\infty\Bigr\}.
\]
Thus (2) is equivalent to $\theta_I=0$. For the $\beta$-L{\"u}roth system, $\theta_{\Phi_I}=\theta_I$. By Theorem~\ref{mu-absolute-regularity}, this is equivalent to absolute regularity of $\Phi_I$, and the conclusion holds.  
\end{proof}

\section{Hausdorff measures of $\beta$-L\"uroth systems}\label{sec_HMD}
In this section,  we study the $h$-dimensional Hausdorff measure of a regular 
$\beta$-L\"uroth system $\Phi_I$ with $\dim_H E_I=h$. 
Recall that $B_k^l=\sum_{n=k}^l B_n$ and $B_k^\infty=\sum_{n=k}^\infty B_n$.

Since $\Phi_I$ is regular  with $\dim_H E_I=h$,  we have  $P(h)=0$ and $\sum_{i\in I}(B_i)^h=1.$
Hence, for every $n\ge1$,
\[
\sum_{\omega\in I^n}\operatorname{diam}(\varphi_\omega(X))^h
=
\Bigl(\sum_{i\in I}(B_i)^h\Bigr)^n
=1,
\]
which implies that $\mathcal H^h(E_I)\le1$. Therefore,   Theorem~\ref{hm} 
characterizes   $\mathcal H^h(E_I)=0$ in terms of  
\[
\sum_{n=k}^{l}\mathbf{1}_I(n)(B_n)^h
\qquad\text{and}\qquad
(B_k^l)^h,
\]
and, equivalently, characterizes the positivity of $\mathcal H^h(E_I)$.

The proof of Theorem~\ref{hm} relies on the following Frostman-type criterion in  \cite{MU99},   which relates the positivity 
of the $h$-dimensional Hausdorff measure to a local upper bound for the 
$h$-conformal measure.

\begin{mainthm}\label{ghm}
Let $\Phi_I=\big\{\varphi_i\big\}_{i\in I}$ be regular with $\dim_H E=h$, and let $m$ be the $h$-conformal probability measure. Then the
following statements are equivalent.
\begin{enumerate}
\item[(1)] $0<\mathcal H^h(E)<\infty$.
\item[(2)] There exist $0<L<\infty$, $\gamma\ge1$, and a finite set
$J\subset I$ such that, for every $i\in I\setminus J$ and every
$r>\gamma B_i$, there exists $x\in\varphi_i(X)$ for which
\[
m(B(x,r))\le Lr^h.
\]
\end{enumerate}
\end{mainthm}

\begin{proof}[Proof of Theorem~\ref{hm}]
Let
\begin{equation}\label{def_M}
M=\sup_{k<l}\frac{1}{(B_k^l)^h}\sum_{n=k}^{l}\mathbf{1}_I(n)(B_n)^h.
\end{equation} 
It suffices to prove the equivalent statement  $0<\mathcal H^h(E_I)<\infty \Longleftrightarrow  M <\infty. $

First, suppose that $0<\mathcal H^h(E_I)<\infty$. By condition~(2) of Theorem~\ref{ghm}, there exist $L>0$, $\gamma\ge1$, and a finite set $J\subset I$ such that, for every $i\in I\setminus J$ and every $r>\gamma B_i$, there exists $x\in\varphi_i(X)$ such that
\[
m(B(x,r))\le Lr^h.
\]
Fix integers $k<l$. We claim that
\[
\sum_{n=k}^{l}\mathbf{1}_I(n)(B_n)^h
\lesssim (B_k^l)^h
\]
with a constant independent of $k$ and $l$. If $I\cap[k,l]=\emptyset$, then
$$
\sum_{n=k}^{l}\mathbf{1}_I(n)(B_n)^h=0<(B_k^l)^h.
$$
 If $I\cap[k,l]\subset J$, since $B_n\le B_k^l$ for every $k\le n\le l$ we have
\[
\sum_{n=k}^{l}\mathbf{1}_I(n)(B_n)^h
=
\sum_{n\in I\cap[k,l]}(B_n)^h
\le
\#J\,(B_k^l)^h.
\]

Suppose that $(I\cap[k,l])\setminus J\not=\emptyset$. Choose
$i\in(I\cap[k,l])\setminus J$ and set
\[
r=(\gamma+1)B_k^l.
\]
Since $B_i\le B_k^l$, we have $r>\gamma B_i$. Thus there exists
$x\in\varphi_i(X)$ such that $m(B(x,r))\le Lr^h$. 
Moreover, the interval $\bigcup_{n=k}^{l}\varphi_n(X)=[B_{l+1}^\infty,B_k^\infty]$ has length \(B_k^l\) and contains \(x\). Since \(r>B_k^l\), it follows that $\bigcup_{n=k}^{l}\varphi_n(X)\subset B(x,r).$
Hence
\[
\bigcup_{n\in I\cap[k,l]}\varphi_n(E_I)\subset\bigcup_{n\in I\cap[k,l]}\varphi_n(X)\subset B(x,r).
\]
By \eqref{conformal measure}, we have $m(\varphi_n(E_I))=(B_n)^h$, and
\begin{equation*}
\begin{aligned}
\sum_{n=k}^{l}\mathbf{1}_I(n)(B_n)^h &=m\Bigl(\bigcup_{n\in I\cap[k,l]}\varphi_n(E_I)\Bigr)\le m(B(x,r))\le L(\gamma+1)^h(B_k^l)^h.
\end{aligned}
\end{equation*}
Therefore $M\leq \max\{\# J, L(\gamma+1)^h\}<\infty$.

Next, suppose  $M<\infty.$ 
Fix $i\in I$,
$x\in\varphi_i(X)$, and $r>B_i$. We claim that there exists 
$C>0$, independent of $i$, $x$, and $r$, such that
\[
m(B(x,r))\le Cr^h.
\]

Set 
$$
S_I =\{j\in I:  \varphi_j(X)\cap B(x,r)\neq\emptyset\}. 
$$
Let $p=\min S_I $. Then $p\le i$. If $p=i$, we have $B_p=B_i<r$. Suppose that $p<i$ and write
$p=i_{j_0}$. 
If $p=i-1$, then $B_{i_{j_0}+1}^{i_{j_0+1}}=B_i<r$.  Otherwise, since $i\in I$, it follows that $i_{{j_0}+1}< i$.
Since the ball intersects $\varphi_p(X)$ and its centre lies in
$\varphi_i(X)$, we have $B_{p+1}^{i-1}<r. $
Moreover, $B_{p+1}^{i_{{j_0}+1}}<2r$ since $B_{p+1}^{i_{{j_0}+1}}\le B_{p+1}^{i-1}<r$ if $i_{{j_0}+1}<i$ and $B_{p+1}^{i_{{j_0}+1}} =B_{p+1}^{i-1}+B_i<2r$ if $i_{{j_0}+1}=i$.
Since $\Phi_I$ satisfies  GRC, by \eqref{hausdorff-hyp-gap}, it follows that
\begin{equation}\label{hausdorff-left-boundary}
B_p
<
2r\sup_{j\ge1}
\frac{B_{i_j}}{B_{i_j+1}^{i_{j+1}}}
\qquad\text{ for every }p<i.
\end{equation}

We then distinguish two cases according to whether $S_I$ is finite or infinite.

(I) Suppose  $\# S_I <\infty$. Let $q=\max S_I $. Then $q\ge i$, and since $\{B_n\}$ is decreasing, we have $B_q\le B_i<r.$ By \eqref{hausdorff-left-boundary} we have
\begin{equation}\label{ineq_Bpqr}
B_p^q \le 2r+B_p+B_q< \biggl( 4+2\sup_{j\ge1} \frac{B_{i_j}}{B_{i_j+1}^{i_{j+1}}} \biggr)r.
\end{equation}
If $p<q$, since $B(x,r)\cap E_I\subset\bigcup_{n=p}^q\varphi_n(E_I),$ we have
\begin{equation*}
m(B(x,r))=m(B(x,r)\cap E_I) \le \sum_{n=p}^{q}\mathbf{1}_I(n)(B_n)^h .
\end{equation*}
Combining it with \eqref{def_M} and \eqref{ineq_Bpqr}, we obtain 
\begin{equation*}
m(B(x,r))\le M(B_p^q)^h \le M\biggl(4+2\sup_{j\ge1}\frac{B_{i_j}}{B_{i_j+1}^{i_{j+1}}} \biggr)^h r^h.
\end{equation*}

If $p=q$, then  we have $p=q=i$ and  
\[
m(B(x,r))\le m(\varphi_i(E))=(B_i)^h<r^h.
\]

(II) Suppose  $\# S_I =\infty$. It is obvious that $0\in\overline{B(x,r)}$. Since $M<\infty$, for every integer $N>p$, we have 
\[
\sum_{n=p}^{N}\mathbf{1}_I(n)(B_n)^h
\le M (B_p^N)^h.
\]
As $N$ tends to infinity, it follows that
\begin{equation}\label{hausdorff-tail-block}
\sum_{n=p}^{\infty}\mathbf{1}_I(n)(B_n)^h \le M(B_p^\infty)^h.
\end{equation}
Since $0\in\overline{B(x,r)}$ and $\varphi_p(X)\cap B(x,r)\neq \emptyset$, we have
\[
B_p^\infty\le2r+B_p.
\]
Note that $B_p<r$ if $p=i$. Combining it with \eqref{hausdorff-left-boundary}, we have that for all $p\leq i$ 
\begin{equation}\label{hausdorff-tail-length}
B_p^\infty
<
\biggl(
3+2\sup_{j\ge1}
\frac{B_{i_j}}{B_{i_j+1}^{i_{j+1}}}
\bigg)r.
\end{equation}
Since $p=\min S_I$ and $B(x,r)\cap E_I \subset \bigcup_{ n\in I, n\ge p}\varphi_n(X),$ it is clear that 
\[
m(B(x,r)) \le \sum_{n\in I,   n\ge p}(B_n)^h.
\]
Combining this with \eqref{hausdorff-tail-block} and 
\eqref{hausdorff-tail-length}, we obtain $m(B(x,r))\lesssim r^h$.

Therefore there exists a constant $C>0$, independent of $i$, $x$, and $r$, such that
\[
m(B(x,r))\le Cr^h
\qquad
\text{for every }i\in I,\ x\in\varphi_i(X),\ r>B_i.
\]
By condition~(2) of Theorem~\ref{ghm}, we have $\mathcal H^h(E_I)>0$. 
\end{proof}

\section{Packing measures of  $\beta$-L{\"u}roth systems}\label{sec:packing}
In this section, we study the packing measure of the attractor and give the proof of  Theorem~\ref{packing} and corollary \ref{cor_PT}. Recall that
\[
\varphi_n(X)=[B_{n+1}^\infty,B_n^\infty], \quad B_k^l=\sum_{n=k}^l B_n, \quad B_k^\infty=\sum_{n=k}^\infty B_n.
\]

The following lemma gives a simple condition for a collection of first-level cylinders with indices from $k$ to $l$ to lie in a ball. 
\begin{lem}\label{packing-block-geometry}
Let $1\leq k\leq i\leq l\leq\infty$ and $r>0$. There exists
$x\in\varphi_i(X)$ such that
\[
\bigcup_{n=k}^{l}\varphi_n(X)\subset B(x,r)
\]
if and only if
\begin{equation}\label{packing-geometry-condition}
B_k^{i-1}<r,
\qquad   
\begin{cases}
B_{i+1}^{l}<r, & l<\infty,\\
B_{i+1}^\infty\leq r, & l=\infty,
\end{cases}       
\qquad
B_k^l<2r.
\end{equation}
\end{lem}

\begin{proof}
We first prove the necessity. Suppose that there exists
$x\in\varphi_i(X)$ such that
\[
\bigcup_{n=k}^{l}\varphi_n(X)\subset B(x,r).
\]
If $l<\infty$, then $\bigcup_{n=k}^{l}\varphi_n(X)=[B_{l+1}^\infty, B_k^\infty]$, so
\[
B_k^\infty-x<r,
\qquad
x-B_{l+1}^\infty<r.
\]
If $l=\infty$, then $\bigcup_{n=k}^{\infty}\varphi_n(X)=(0,B_k^\infty]$, so
\[
B_k^\infty-x<r,
\qquad
x\leq r.
\]
Since $x\leq B_i^\infty$ and $x\geq B_{i+1}^\infty$, it follows that
$B_k^{i-1}=B_k^\infty-B_i^\infty<r$, and that
$B_{i+1}^{l}<r$ if $l<\infty$, while $B_{i+1}^\infty\leq r$ if $l=\infty$.
Moreover, $B_k^l<2r$ follows from the two endpoint bounds in each case. Thus
\eqref{packing-geometry-condition} holds.

For sufficiency, suppose that \eqref{packing-geometry-condition} holds.
If $l<\infty$, then
\[
\bigcup_{n=k}^{l}\varphi_n(X)=[B_{l+1}^\infty,B_k^\infty].
\]
By \eqref{packing-geometry-condition}, we have
\[
B_k^\infty-r<B_i^\infty,
\qquad
B_{i+1}^\infty<B_{l+1}^\infty+r,
\qquad
B_k^\infty-r<B_{l+1}^\infty+r.
\]
Hence we can choose
\[
x\in[B_{i+1}^\infty,B_i^\infty]
\cap(B_k^\infty-r,\;B_{l+1}^\infty+r),
\]
which gives $\bigcup_{n=k}^{l}\varphi_n(X)\subset B(x,r)$.

If $l=\infty$, then $\bigcup_{n=k}^{\infty}\varphi_n(X)=(0,B_k^\infty]$.
By \eqref{packing-geometry-condition}, we have
\[
B_k^\infty-r<B_i^\infty,
\qquad
B_{i+1}^\infty\leq r,
\qquad
B_k^\infty-r<r.
\]
Therefore, we can choose
\[
x\in[B_{i+1}^\infty,B_i^\infty]\cap(B_k^\infty-r,\;r].
\]
This $x$ satisfies $\bigcup_{n=k}^{\infty}\varphi_n(X)\subset B(x,r)$.
\end{proof}

In contrast to the Hausdorff measure, where a local upper bound for the 
conformal measure yields positivity, the packing measure is governed by a 
local lower bound. The proof of Theorem~\ref{packing} relies on the 
following dual Frostman-type criterion in \cite{MU99}.
\begin{mainthm}\label{pm}
Let $\Phi_I=\{\varphi_i\}_{i\in I}$ be regular, let $\dim_H E=h$, and let
$m$ be the $h$-conformal probability measure. Then the following statements
are equivalent.
\begin{enumerate}
\item[(1)] $\mathcal P^h(E)<\infty$.
\item[(2)] There exist $L>0$ and $\xi>0$ such that for every $\gamma\geq1$,
there exists a finite set $J\subset I$ with the following property: for every
$i\in I\setminus J$ and every $r$ satisfying
\[
\gamma B_i<r\leq\xi,
\]
there exists $x\in\varphi_i(X)$ such that
\[
m(B(x,r))\geq Lr^h.
\]
\end{enumerate}
\end{mainthm}

\begin{proof}[Proof of Theorem~\ref{packing}]
Let $m$ be the $h$-conformal measure. Since the system is affine, we have
\begin{equation}\label{packing-cylinder-mass}
m(\mathbf{B}_n)=(B_n)^h,
\qquad n\in I.
\end{equation}

First, we prove the sufficiency. Fix
$i\in I\setminus J$ and $B_i<r\le\xi$, and choose $k,l$ satisfying
\eqref{packing-main-geometry} and \eqref{packing-main-mass}. By
Lemma~\ref{packing-block-geometry}, there exists $x\in\varphi_i(X)$ such that
\[
\bigcup_{n=k}^{l}\varphi_n(X)\subset B(x,r).
\]
Therefore, by \eqref{packing-cylinder-mass}, we have
\begin{equation*}
\begin{aligned}
m(B(x,r)) \geq \sum_{n=k}^{l}\mathbf{1}_I(n)m(\varphi_n(E_I))=\sum_{n=k}^{l}\mathbf{1}_I(n)(B_n)^h\geq \kappa r^h.
\end{aligned}
\end{equation*}
Thus condition~(2) of Theorem~\ref{pm} holds for $\gamma=1$. 
Since $\gamma\ge1$, we have $\gamma B_i\ge B_i$, and therefore
\[
(\gamma B_i,\xi]\subseteq (B_i,\xi].
\]
Thus, if $i\in I\setminus J$ and $r$ satisfies $\gamma B_i<r\le\xi$, 
then automatically $B_i<r\le\xi$. By the case $\gamma=1$ already proved, 
there exists $x\in\varphi_i(X)$ such that
\[
m(B(x,r))\ge \kappa r^h.
\]
Hence condition~(2) of Theorem~\ref{pm} holds for this $\gamma$, with the 
same constants $\kappa,\xi$ and the same finite set $J$. Since $\gamma\ge1$ 
was arbitrary, by Theorem~\ref{pm}, we obtain 
\[
\mathcal P^h(E_I)<\infty.
\]

For necessity, suppose that
$\mathcal P^h(E_I)<\infty$. By Theorem~\ref{pm} with $\gamma=1$, there exist $L>0$, $\xi>0$, and a finite set $J\subset I$ such that, for every $i\in I\setminus J$ and every $B_i<r\le\xi$, there exists $x\in\varphi_i(X)$ satisfying
\begin{equation}\label{packing-mu-lower}
m(B(x,r))\geq Lr^h.
\end{equation}
By \eqref{hausdorff-hyp-gap}, choose $K\ge1$ such that
\begin{equation}\label{packing-gap-K}
B_{i_j}\leq K B_{i_j+1}^{i_{j+1}}
\qquad\text{for every }j\geq1.
\end{equation}

Fix $i\in I\setminus J$ and $B_i<r\le\xi$, and choose $x$ as in
\eqref{packing-mu-lower}.  Set
$$
S_I=\{j\in I:  \varphi_j(X)\cap B(x,r)\neq\emptyset\}. 
$$
Since $x\in\varphi_i(X)\cap B(x,r)$, we have $i\in S_I$ and 
$S_I \neq\emptyset$. Let $k=\min S_I $, and let $l=\max S_I$ 
if $S_I$ is finite and $l=\infty$ otherwise. Then 
we have 
\[
B(x,r)\cap E_I\subset\bigcup_{n=k}^{l}\varphi_n(E_I).
\]
Therefore, by \eqref{packing-cylinder-mass}, it follows that
\begin{equation}\label{packing-ball-block-upper}
\begin{aligned}
m(B(x,r)) &\leq \sum_{n=k}^{l}\mathbf{1}_I(n)m(\mathbf{B}_n)
=\sum_{n=k}^{l}\mathbf{1}_I(n)(B_n)^h.
\end{aligned}
\end{equation}
Combining \eqref{packing-mu-lower} and
\eqref{packing-ball-block-upper}, we obtain
\begin{equation}\label{packing-block-mass-at-r}
\sum_{n=k}^{l}\mathbf{1}_I(n)(B_n)^h\geq Lr^h.
\end{equation}

We claim that
\[
\bigcup_{n=k}^{l}\varphi_n(X) \subset B\bigl(x,(2K+1)r\bigr).
\]
Since $\bigcup_{n=k}^{l}\varphi_n(X)\subset [B_{l+1}^\infty,B_k^\infty],$
with the convention $B_\infty^\infty:=0$, it suffices to prove  
\[
\text{(I) } B_k^\infty-x<(2K+1)r,
\qquad
\text{(II) } x-B_{l+1}^\infty<(2K+1)r.
\]

(I) If $k=i$, then $x\in\varphi_i(X)=[B_{i+1}^\infty,B_i^\infty]$, and hence
\[
B_k^\infty-x =B_i^\infty-x  \le B_i^\infty-B_{i+1}^\infty
=B_i <r <(2K+1)r.
\]

For $k<i$, and write $k=i_j$. Since   $\varphi_k(X)\cap B(x,r)\neq\emptyset$,  it is clear that
\begin{equation}\label{packing-k-left}
B_{k+1}^\infty-x<r.
\end{equation}
Also, since $x\in\varphi_i(X)$,
\begin{equation}\label{packing-x-to-i}
x-B_{i+1}^\infty\le B_i<r.
\end{equation}
Since $i_{j+1}\le i$, we have $B_{k+1}^{i_{j+1}}\le B_{k+1}^{i}.$ Therefore, by \eqref{packing-k-left} and \eqref{packing-x-to-i},
\[
B_{k+1}^{i_{j+1}}
\le B_{k+1}^{i}
=B_{k+1}^\infty-B_{i+1}^\infty
=\bigl(B_{k+1}^\infty-x\bigr)+\bigl(x-B_{i+1}^\infty\bigr)
<r+r=2r.
\]
By \eqref{packing-gap-K}, it follows that $B_k\le K B_{k+1}^{i_{j+1}}<2Kr.$
Combining this with \eqref{packing-k-left}, we obtain
\[
B_k^\infty-x
=B_k+\bigl(B_{k+1}^\infty-x\bigr)
<2Kr+r=(2K+1)r. 
\]
 
(II) If $l<\infty$, then $l\ge i$ since $i\in S_I$ and $l=\max S_I $.  
Since $\{B_n\}$ is decreasing,  it is clear that $B_l\le B_i<r.$ Choose $y\in\varphi_l(X)\cap B(x,r)$, and we have 
\[
x-B_{l+1}^\infty
=(x-y)+\bigl(y-B_{l+1}^\infty\bigr)
<r+B_l<2r<(2K+1)r.
\]
For $l=\infty$, since  $0\in\overline{B(x,r)}$, it immediately follows that
\[
x-B_\infty^\infty=x\le r<(2K+1)r.
\]

  Applying Lemma~\ref{packing-block-geometry} with radius $(2K+1)r$,
and noting that $B_{i+1}^\infty\le x\le r<(2K+1)r$ when $l=\infty$,
we obtain
\begin{equation}\label{packing-enlarged-geometry}
B_k^{i-1}<(2K+1)r, \qquad B_{i+1}^{l}<(2K+1)r, \qquad B_k^l<2(2K+1)r.
\end{equation}

Finally, we have proved that, for every $i\in I\setminus J$ and every $u$ 
with $B_i<u\le\xi$, there exist $k\le i\le l$ such that
\begin{equation}\label{packing-auxiliary-conclusion}
B_k^{i-1}<(2K+1)u, \qquad B_{i+1}^{l}<(2K+1)u, \qquad B_k^l<2(2K+1)u,
\end{equation}
and
\begin{equation}\label{packing-auxiliary-mass}
\sum_{n=k}^{l}\mathbf{1}_I(n)(B_n)^h\ge Lu^h.
\end{equation}
It remains to prove the same conclusion for every $r$ with 
$B_i<r\le\xi$, with  
\[
B_k^{i-1}<r,\qquad B_{i+1}^{l}<r, \qquad B_k^l<2r.
\]

 Fix $i\in I\setminus J$ and $B_i<r\le\xi$. If $r>(2K+1)B_i$, set $u=\frac{r}{2K+1}$, and we have  $B_i<u\le\xi.$  By \eqref{packing-auxiliary-conclusion} and \eqref{packing-auxiliary-mass}, 
there exist $k\le i\le l$ such that
\[
B_k^{i-1}<(2K+1)u=r,
\qquad
B_{i+1}^{l}<(2K+1)u=r,
\qquad
B_k^l<2(2K+1)u=2r,
\]
and
\[
\sum_{n=k}^{l}\mathbf{1}_I(n)(B_n)^h
\ge Lu^h
=\frac{L}{(2K+1)^h}r^h.
\]
If $B_i<r\le(2K+1)B_i$, take $k=l=i$, and we have 
\[
B_k^{i-1}=B_{i+1}^{l}=0<r,
\qquad
B_k^l=B_i<r<2r.
\]
Moreover,  since $r\le(2K+1)B_i$, we have $B_i\ge r/(2K+1)$,  and hence
\[
\sum_{n=k}^{l}\mathbf{1}_I(n)(B_n)^h
=(B_i)^h
\ge (2K+1)^{-h}r^h.
\]

Combining the two cases, we conclude that for every $i\in I\setminus J$ 
and every $r$ with $B_i<r\le\xi$, there exist $k\le i\le l$ such that
\[
B_k^{i-1}<r,
\qquad
B_{i+1}^{l}<r,
\qquad
B_k^l<2r,
\]
and
\[
\sum_{n=k}^{l}\mathbf{1}_I(n)(B_n)^h
\ge \kappa r^h,
\]
where $\kappa=\frac{\min\{L,1\}}{(2K+1)^h}.$  
\end{proof}

\begin{proof}[Proof of Corollary~\ref{cor_PT}]
For each integer $a\geq 1$, we write 
\[
f(a)=\frac{1}{(B_a^\infty)^h}\sum_{n=a}^{\infty}\mathbf{1}_I(n)(B_n)^h.
\]
Since $\mathcal P^h(E_I)<\infty$, let  $\kappa>0$, $\xi>0$, and a finite set $J\subset I$  be chosen as in  Theorem~\ref{packing}. 

Choose $a$ sufficiently large that $\frac12 B_a^\infty\le\xi.$
Since $I$ is infinite and $\lim_{n\to\infty}B_n=0$, we choose 
$i\in I\setminus J$ with $i>a$ so large that
\begin{equation}\label{packing-tail-choice-i}
B_i<\frac12 B_a^\infty
\qquad\text{and}\qquad
B_a^{i-1}>\frac12 B_a^\infty.
\end{equation}
Set $r=\frac12 B_a^\infty$. Then $B_i<r\le\xi$, so 
\eqref{packing-main-geometry}--\eqref{packing-main-mass} apply and give 
indices $k\le i\le l$ such that $B_k^{i-1}<r$ and
\[
\sum_{n=k}^{l}\mathbf{1}_I(n)(B_n)^h
\ge \kappa r^h.
\]

We claim that $k>a$. Indeed, if $k\le a$, then by 
\eqref{packing-tail-choice-i}, $B_k^{i-1}\ge B_a^{i-1}>r,$
contradicting $B_k^{i-1}<r$. Hence $k>a$. Consequently 
$[k,l]\subset[a,\infty)$, and therefore
\[
\sum_{n=a}^{\infty}\mathbf{1}_I(n)(B_n)^h
\ge
\sum_{n=k}^{l}\mathbf{1}_I(n)(B_n)^h
\ge \kappa r^h
= \kappa 2^{-h}(B_a^\infty)^h.
\]
Thus, for all sufficiently large $a$, we have 
\begin{equation}\label{packing-tail-lower-large-a}
f(a) \ge \kappa 2^{-h}. 
\end{equation}

Let $a_0$ be such that \eqref{packing-tail-lower-large-a} holds for every 
$a\ge a_0$. Assume  $a_0>1$. Then for $1\le a<a_0$, the quantity $f(a)$ 
is finite and positive. Indeed,
\[
\sum_{n=a}^{\infty}\mathbf{1}_I(n)(B_n)^h
\ge \sum_{n=a_0}^{\infty}\mathbf{1}_I(n)(B_n)^h,
\]
and $B_a^\infty\le1$. Hence
\[
\inf_{a\ge1}f(a)\ge \min\{\sum_{n=a_0}^{\infty}\mathbf{1}_I(n)(B_n)^h,\kappa 2^{-h}\}>0.
\]
If $a_0=1$, then \eqref{packing-tail-lower-large-a} already holds for all  $a\geq 1$, and the conclusion follows.
\end{proof}

\section{Tail-regularity of  $\beta$-L{\"u}roth systems}\label{sec:tail-regular}
The criteria in Theorems~\ref{hm} and~\ref{packing} are stated in terms of 
the block sums $(B_k^l)^h$ and $\sum_{n=k}^{l}\mathbf{1}_I(n)(B_n)^h$, 
and are therefore difficult to verify directly for a given index set $I$. 
However, if the partition $\beta$ is tail-regular, these block sums are 
governed by a single arithmetic scale $R_n=\frac{B_n^\infty}{B_n}.$ Under this assumption, the criteria of Theorems~\ref{hm} and~\ref{packing} 
reduce to counting conditions for $I$ on intervals of length comparable 
to $R_n$, as stated in Theorems~\ref{tail-hausdorff} and~\ref{tail-packing}; 
the proofs are given in Subsection~\ref{subsec_HPM}. In the final subsection, 
we verify the tail-regularity condition for regularly varying and 
stretched-exponential partitions.

\subsection{The tail scale and its consequences}
The following lemma collects the two basic consequences of tail-regularity 
used throughout this section: a uniform bound for the gap ratios 
$B_{i_j}/B_{i_j+1}^{i_{j+1}}$, and a geometric decay of the tails 
$B_n^\infty$ at the natural scale $L_n=\lfloor\delta R_n\rfloor$.
\begin{lem}\label{tail-regular-lemma}
Assume that the partition $\beta$ is tail-regular. Then the following   hold.
\begin{enumerate}
\item[(a)] For every infinite $I$,
\[
\sup_{j\ge1}\frac{B_{i_j}}{B_{i_j+1}^{i_{j+1}}}<\infty.
\]
\item[(b)] There exist $\delta<\eta$, $q\in(0,1)$, and $n_1\in\mathbb N$ such that, for $L_n=\lfloor\delta R_n\rfloor, n\ge n_1,$
one has
\begin{equation}\label{tail-natural-block}
B_n^{n+L_n-1}\asymp B_n^\infty
\end{equation}
and
\begin{equation}\label{tail-geometric-decay}
B_{n+L_n}^\infty\le qB_n^\infty.
\end{equation}
\end{enumerate}
\end{lem}

\begin{proof}
(a). Recall that $R_n=\frac{B_n^\infty}{B_n}.$ Since $\beta$ is tail-regular, we have $\lim_{n\to\infty}R_n=\infty$ and there exist constants $A\ge1$, $\eta>0$, and $n_0\in\mathbb N$ such that, for every $n\ge n_0$, $\eta R_n\le\frac n2$,
and
\begin{equation}\label{tail-local-comparison-lemma}
A^{-1}B_n\le B_m\le AB_n
\end{equation}
whenever $m\in\mathbb N$ and $|m-n|\le\eta R_n$.

 Since $\lim_{n\to\infty}R_n=\infty$, there exists $N\ge n_0$ such that $1\le\eta R_n$ for every $n\ge N$. By \eqref{tail-local-comparison-lemma}, we have  $B_{n+1}\ge A^{-1}B_n$  for all $n\ge N$.
Let $i_j\ge N$. Since $B_{i_j+1}^{i_{j+1}}\ge B_{i_j+1}$, we have
\[
\frac{B_{i_j}}{B_{i_j+1}^{i_{j+1}}}
\le
\frac{B_{i_j}}{B_{i_j+1}}
\le A.
\]
The remaining ratios correspond to only finitely many $i_j<N$, and each of them is finite. Therefore conclusion (a) follows.    

(b). Choose $\delta>0$ such that $0<\delta<\min\Bigl\{\eta,\frac{1}{4A}\Bigr\}$. Since $\lim_{n\to\infty}R_n=\infty$, there exists $n_1\ge n_0$ such that
\[
L_n=\lfloor\delta R_n\rfloor\ge\frac{\delta}{2}R_n
\qquad(n\ge n_1).
\]
For $0\le j<L_n$, we have $j\le\delta R_n\le\eta R_n$. Hence, by \eqref{tail-local-comparison-lemma},
\[
A^{-1}B_n\le B_{n+j}\le AB_n.
\]
It follows that
\begin{equation}\label{tail-natural-block-bounds}
\frac{\delta}{2A}R_nB_n
\le
B_n^{n+L_n-1}
\le
A\delta R_nB_n.
\end{equation}
Recalling that $R_n=\frac{B_n^\infty}{B_n}$, we have $R_nB_n=B_n^\infty$. Thus \eqref{tail-natural-block} follows.

Finally,
\[
B_{n+L_n}^\infty=B_n^\infty-B_n^{n+L_n-1}.
\]
Using the lower bound in \eqref{tail-natural-block-bounds}, we obtain
\[
B_{n+L_n}^\infty
\le
\Bigl(1-\frac{\delta}{2A}\Bigr)B_n^\infty.
\]
Taking $q=1-\frac{\delta}{2A}$ gives \eqref{tail-geometric-decay}, and the proof is complete.
\end{proof}

\subsection{Hausdorff and packing measures under tail-regularity} \label{subsec_HPM}
In this subsection we prove Theorems~\ref{tail-hausdorff} and~\ref{tail-packing}, 
Throughout, we assume that the partition $\beta$ is tail-regular,   and that $\Phi_I$ is a regular $\beta$-L\"uroth system with 
$\dim_H E_I=h$.

First, we give the proof of Theorem~\ref{tail-hausdorff}, which characterizes the positivity of $\mathcal H^h(E_I)$ by an upper counting condition for $I$ on intervals of length up to the tail scale $\delta R_n$.

\begin{proof}[Proof of Theorem~\ref{tail-hausdorff}]
Since $\beta$ is tail-regular, we have $\lim_{n\to\infty}R_n=\infty$ and there exist constants $A\ge1$, $\eta>0$, and $n_0\in\mathbb N$ such that, for every $n\ge n_0$, $\eta R_n\le\frac n2$, and
\begin{equation}\label{tail-local-comparison-H}
A^{-1}B_n\le B_m\le AB_n
\end{equation}
whenever $m\in\mathbb N$ and $|m-n|\le\eta R_n$.
Let $\delta>0$ be chosen as in Lemma~\ref{tail-regular-lemma}. By Lemma~\ref{tail-regular-lemma}(a), condition \eqref{hausdorff-hyp-gap} holds, so Theorem~\ref{hm} applies.

Suppose $\mathcal H^h(E_I)>0$. By Theorem~\ref{hm}, there exists $M>0$ 
such that for all $k<l$,
\begin{equation}\label{tail-H-block-bound}
\sum_{j=k}^{l}\mathbf{1}_I(j)(B_j)^h
\le
M(B_k^l)^h.
\end{equation}
We claim that there exists $C>0$ such that, for every sufficiently large 
$n$ and every integer $1\le m\le\delta R_n$,
\begin{equation}\label{tail-H-counting-claim}
\#\bigl(I\cap[n,n+m-1]\bigr)\le Cm^h.
\end{equation}

Fix such an $n$ and $m$. If $m=1$, then 
$\#(I\cap[n,n])\le1$, so \eqref{tail-H-counting-claim} holds with any 
$C\ge1$. Suppose that $m\ge2$. Since $m\le\delta R_n$ and 
$\delta<\eta$, we have
\[
|j-n|\le m-1<\eta R_n
\qquad(n\le j\le n+m-1).
\]
Thus \eqref{tail-local-comparison-H} gives
\[
B_j\ge A^{-1}B_n
\qquad(n\le j\le n+m-1),
\]
and hence
\begin{equation}\label{tail-H-local-numerator}
\sum_{j=n}^{n+m-1}\mathbf 1_I(j)(B_j)^h
\ge
A^{-h}(B_n)^h\#\bigl(I\cap[n,n+m-1]\bigr).
\end{equation}
On the other hand,
\begin{equation}\label{tail-H-local-denominator}
B_n^{n+m-1}
\le
AmB_n.
\end{equation}
Applying \eqref{tail-H-block-bound} with $k=n$ and $l=n+m-1$, and using 
\eqref{tail-H-local-numerator} and \eqref{tail-H-local-denominator}, we obtain
\[
A^{-h}(B_n)^h\#\bigl(I\cap[n,n+m-1]\bigr)
\le
M(AmB_n)^h
=
MA^h m^h B_n^h.
\]
Dividing by $A^{-h}B_n^h>0$, we get
\[
\#\bigl(I\cap[n,n+m-1]\bigr)
\le
MA^{2h}m^h.
\]
By setting $C=\max\{MA^{2h},1\}$,   the claim holds for all $1\le m\le\delta R_n$.

 Suppose that \eqref{tail-H-counting-intro} holds. By Theorem~\ref{hm}, it is sufficient to prove 
\begin{equation}\label{tail-H-global-block}
\sup_{k<l}
\frac{ 1}
{(B_k^l)^h} \sum_{n=k}^{l}\mathbf{1}_I(n)(B_n)^h<\infty.
\end{equation}
Fix $k<l$ with $k$ sufficiently large. Set $m=l-k+1$.

First, suppose that $m\le\delta R_k$. By \eqref{tail-H-counting-intro} and \eqref{tail-local-comparison-H}, we have 
\begin{equation*}
\sum_{n=k}^{l}\mathbf{1}_I(n)(B_n)^h \le CA^h m^h(B_k)^h.
\end{equation*}
Combining it with  $B_k^l\ge A^{-1}mB_k$,  we obtain \eqref{tail-H-global-block}.

Suppose  that $m>\delta R_k$. Set $n_0=k$, and for each $j\ge0$, 
define
\[
n_{j+1}=n_j+L_{n_j},
\qquad
L_n=\lfloor\delta R_n\rfloor.
\]
We first bound the tail mass. By \eqref{tail-H-counting-intro} and 
\eqref{tail-local-comparison-H},
\[
\sum_{\substack{n_j\le n<n_{j+1}\\n\in I}}(B_n)^h
\le
CA^hL_{n_j}^h(B_{n_j})^h.
\]
Since $L_n\le (\delta+1)R_n$, we have $L_{n_j}B_{n_j}\le(\delta+1)B_{n_j}^\infty$. Hence
\[
CA^hL_{n_j}^h(B_{n_j})^h
\lesssim
(B_{n_j}^\infty)^h.
\]
Moreover, Lemma~\ref{tail-regular-lemma}(b) gives 
$B_{n_{j+1}}^\infty\le qB_{n_j}^\infty$ for some $q\in(0,1)$. Since 
$\sum_{n\in I}(B_n)^h=1$ and $I$ is infinite, we have $h>0$. Therefore
\begin{equation}\label{tail-H-tail-mass}
\begin{aligned}
\sum_{\substack{n\ge k\\n\in I}}(B_n)^h
&=
\sum_{j=0}^{\infty}
\sum_{\substack{n_j\le n<n_{j+1}\\n\in I}}(B_n)^h
\lesssim
\sum_{j=0}^{\infty}(B_{n_j}^\infty)^h
\le
(B_k^\infty)^h\sum_{j=0}^{\infty}q^{jh}
\lesssim
(B_k^\infty)^h.
\end{aligned}
\end{equation}

Since $m>\delta R_k$, we have $L_k=\lfloor\delta R_k\rfloor<m$, and 
$[k,k+L_k-1]\subset[k,l]$. By \eqref{tail-natural-block},
\[
B_k^l\ge B_k^{k+L_k-1}\gtrsim B_k^\infty.
\]
Combining this with \eqref{tail-H-tail-mass}, we obtain
\[
\frac{\sum_{n=k}^{l}\mathbf{1}_I(n)(B_n)^h}{(B_k^l)^h}
\lesssim 1.
\]
 
 For the finitely many remaining $k$, regularity gives 
$\sum_{n\in I}(B_n)^h=1$, and $B_k^l\ge B_k>0$. Hence
\[
\frac{\sum_{n=k}^{l}\mathbf{1}_I(n)(B_n)^h}{(B_k^l)^h}
\le
\frac{1}{(B_k)^h}<\infty.
\]

Therefore, it is uniformly bounded for all $k$ as well, and \eqref{tail-H-global-block} follows.
\end{proof}

Next, we give the proof of Theorem~\ref{tail-packing}, where the criterion consists of a lower counting 
condition for $I$ and a tail condition for the mass of $I$ beyond $n$.

\begin{proof}[Proof of Theorem~\ref{tail-packing}]
Since $\beta$ is tail-regular, we have $\lim_{n\to\infty}R_n=\infty$ and there exist constants $A\ge1$, $\eta>0$, and $n_0\in\mathbb N$ such that, for every $n\ge n_0$, $\eta R_n\le\frac n2$, and
\begin{equation}\label{tail-local-comparison-P}
A^{-1}B_n\le B_m\le AB_n,
\end{equation}
whenever $m\in\mathbb N$ and $|m-n|\le\eta R_n$.
Let $\delta>0$ be chosen as in Lemma~\ref{tail-regular-lemma}. By Lemma~\ref{tail-regular-lemma}(a), condition \eqref{hausdorff-hyp-gap} holds, so Theorem~\ref{packing} applies.

First, suppose that $\mathcal P^h(E_I)<\infty$. Condition (ii) follows directly from Corollary~\ref{cor_PT}. It remains to prove condition (i).

Let $\kappa>0$, $\xi>0$, and the finite set $J\subset I$ be given by Theorem~\ref{packing}. Fix a sufficiently large $i\in I\setminus J$ and an integer $1\le m\le\delta R_i$. We consider  two cases: (a) $m>4A$; (b) $1\leq m\leq 4A.$ 

(a)  Set $r=\frac{mB_i}{2A}.$  Then $r>B_i$. Since $m\le\delta R_i$ and $\lim_{i\to\infty}B_i^\infty=0$, we   have $r\le\xi$ for all sufficiently large $i$. By Theorem~\ref{packing}, there exist $k\le i\le l$ satisfying \eqref{packing-main-geometry} and \eqref{packing-main-mass}.

We claim that
\begin{equation*}
i-m<k\le i\le l<i+m.
\end{equation*}
Suppose that $k\le i-m$. Since $m\le\delta R_i\le\eta R_i$, \eqref{tail-local-comparison-P} gives
\[
B_k^{i-1}
\ge
B_{i-m}^{i-1}
\ge
\frac{m}{A}B_i
>r,
\]
contradicting $B_k^{i-1}<r$ from  \eqref{packing-main-geometry}. Hence $k>i-m$. By the same argument, we have  \[
B_{i+1}^{i+m}\ge\frac{m}{A}B_i>r,
\]
and it follows that  $l<i+m$. Hence, the claim follows.

By \eqref{tail-local-comparison-P}, we have $B_n\le AB_i$ for every $k\le n\le l$. 
By  \eqref{packing-main-mass}, it follows that 
\[
A^hB_i^h\#\bigl(I\cap[i-m,i+m]\bigr)
\ge
\kappa r^h.
\]
Therefore
\begin{equation}\label{tail-P-local-lower-large-m}
\#\bigl(I\cap[i-m,i+m]\bigr)
\ge
\frac{\kappa}{2^hA^{2h}}m^h.
\end{equation}

(b) For $1\le m\le4A$.
Since $i\in I$, we have
\[
\#\bigl(I\cap[i-m,i+m]\bigr)\ge1\ge(4A)^{-h}m^h.
\]

 Combining the two cases, we obtain condition~(i) with
\[
c=\min\Bigl\{\frac{\kappa}{2^hA^{2h}},\,(4A)^{-h}\Bigr\}>0.
\]

Next, suppose that conditions (i) and (ii) hold. We verify the condition in Theorem~\ref{packing}. Let $c,c'>0$ be the constants in (i) and (ii). Choose $N$ so large that the tail-regular estimates and both arithmetic conditions hold for every index at least $N$. Increase $N$ further so that
\begin{equation}\label{tail-P-ratio-half}
R_n\ge2,
\qquad
\delta R_n>2
\qquad(n\ge N).
\end{equation}
Let $F=I\cap[1,N+1]$ and choose $0<\xi<B_{N+1}$.

Fix $i\in I\setminus F$ and $B_i<r\le\xi$. Since  $0<\delta<\min \{\eta,\frac{1}{4A}\}$, by \eqref{tail-P-ratio-half}, we have  
\[
8AB_i<4A\delta B_i^\infty<B_i^\infty.
\]
Hence the following four cases cover all possible values of $r$.

First, suppose that $B_i<r\le8AB_i.$ Take $k=l=i$. Then \eqref{packing-main-geometry} holds, and
\[
\sum_{n=k}^{l}\mathbf{1}_I(n)(B_n)^h
=B_i^h
\ge
(8A)^{-h}r^h.
\]

Second, suppose that $8AB_i<r\le4A\delta B_i^\infty.$ Set
\[
m=\Bigl\lfloor\frac{r}{4AB_i}\Bigr\rfloor,
\qquad
k=i-m,
\qquad
l=i+m.
\]
Since $B_i^\infty=R_iB_i$, it is clear that $\delta R_i\ge m\ge\frac{r}{8AB_i}.$
By \eqref{tail-local-comparison-P}, it follows that
\[
B_k^{i-1}\le AmB_i\le\frac r4,
\qquad
B_{i+1}^{l}\le AmB_i\le\frac r4.
\]
Moreover,
\[
B_k^l
\le
B_i+2AmB_i
\le
B_i+\frac r2
<
\frac32r<2r.
\]
Thus \eqref{packing-main-geometry} holds. Since $\delta R_i\ge m\ge\frac{r}{8AB_i}$, by condition (i) and \eqref{tail-local-comparison-P}, we have that
\begin{equation*}
\begin{aligned}
\sum_{n=k}^{l}\mathbf{1}_I(n)(B_n)^h
&\ge
A^{-h}B_i^h\#\bigl(I\cap[i-m,i+m]\bigr)\\
&\ge
cA^{-h}m^hB_i^h\\
&\ge
c(8A^2)^{-h}r^h.
\end{aligned}
\end{equation*}

Third, suppose that $4A\delta B_i^\infty<r\le B_i^\infty.$ Set $m=L_i=\lfloor\delta R_i\rfloor$ and $k=i-m$, $l=i+m$. Choosing sufficiently large $N$, we have $m\le\delta R_i\le\eta R_i$ and $m\ge\frac{\delta}{2}R_i$ for all $i>N$.  By \eqref{tail-local-comparison-P},
\[
B_k^{i-1}\le AmB_i\le A\delta B_i^\infty<\frac r4,
\qquad
B_{i+1}^{l}\le A\delta B_i^\infty<\frac r4.
\]
Also,
\[
B_k^l\le B_i+2A\delta B_i^\infty<\frac32r<2r.
\]
By condition (i),
$$
\sum_{n=k}^{l}\mathbf{1}_I(n)(B_n)^h\ge cA^{-h}m^hB_i^h\ge cA^{-h}\Bigl(\frac{\delta}{2}\Bigr)^h(B_i^\infty)^h\ge cA^{-h}\Bigl(\frac{\delta}{2}\Bigr)^h r^h.
$$

Finally, suppose that $r>B_i^\infty.$
Choose $k\le i$ as the smallest integer satisfying
\begin{equation}\label{tail-P-k-choice}
B_k^\infty<r+B_i^\infty.
\end{equation}
Since $i\ge N+2$ and $r\le\xi<B_{N+1}$, we have 
\[
r+B_i^\infty
\le r+B_{N+2}^\infty
<B_{N+1}^\infty.
\]
Hence $k\ge N+2$, and we have  $B_{k-1}^\infty\ge r+B_i^\infty.$
Using \eqref{tail-P-ratio-half},
\begin{equation}\label{tail-P-Tk-lower}
B_k^\infty
=B_{k-1}^\infty-B_{k-1}
=B_{k-1}^\infty\Bigl(1-\frac1{R_{k-1}}\Bigr)
\ge\frac12B_{k-1}^\infty
\ge\frac r2.
\end{equation}
Take $l=\infty$. Since $r>B_i^\infty$, by  \eqref{tail-P-k-choice}, we obtain   
\[
B_k^{i-1}=B_k^\infty-B_i^\infty<r,
\qquad
B_{i+1}^{\infty}=B_{i+1}^\infty<B_i^\infty<r,
\qquad
B_k^\infty<r+B_i^\infty<2r.
\]
Thus \eqref{packing-main-geometry} holds. By condition (ii) and \eqref{tail-P-Tk-lower},
\begin{equation*}
\sum_{\substack{n\in I\\n\ge k}}(B_n)^h
\ge
c'(B_k^\infty)^h
\ge
c'2^{-h}r^h.
\end{equation*}

Set $\kappa=\min \{(8A)^{-h},  c(8A^2)^{-h}, cA^{-h} (\frac{\delta}{2} )^h, c'2^{-h} \}>0.$      
Combining the four cases, we conclude that for every $i\in I\setminus F$ 
and every $r$ with $B_i<r\le\xi$, there exist $k\le i\le l$ such that 
\eqref{packing-main-geometry} holds and
\[
\sum_{n=k}^{l}\mathbf{1}_I(n)(B_n)^h\ge \kappa r^h.
\]
By Theorem~\ref{packing}, it follows that  $\mathcal P^h(E_I)<\infty.$
 \end{proof}

\subsection{Regularly varying and stretched-exponential partitions}
We say  a positive function $L$ is slowly varying if
\[
\lim_{x\to\infty}\frac{L(cx)}{L(x)}=1
\]
for every $c>0$. We write $a_n\sim b_n$ if
\[
\lim_{n\to\infty}\frac{a_n}{b_n}=1.
\]
 
\begin{prop}\label{regularly-varying-application}
Suppose that there exists a slowly varying function $L$ such that 
\begin{equation}\label{regularly-varying-B}
B_n\sim c n^{-p}L(n),
\end{equation}
for some $p>1$, $c>0$.
Then the partition $\beta$ is tail-regular and
\begin{equation}\label{regularly-varying-R}
R_n\sim\frac{n}{p-1}.
\end{equation}
\end{prop}

\begin{proof}
By Karamata's theorem \cite{BGT87}, it is clear that 
\[
B_n^\infty=\sum_{j=n}^{\infty}B_j
\sim
\frac{c}{p-1}n^{1-p}L(n).
\]
Together with \eqref{regularly-varying-B}, this gives
\[
R_n=\frac{B_n^\infty}{B_n}\sim\frac{n}{p-1},
\]
which proves \eqref{regularly-varying-R}. In particular, $\lim_{n\to\infty}R_n=\infty$ and $R_n\asymp n$. Fix  sufficiently small $\eta>0$. Then \eqref{regularly-varying-R} gives $\eta R_n\le\frac n2$ for all sufficiently large $n$.

To prove $\beta$ is tail-regular, it remains to verify $B_m\asymp  B_n$  for all $m\in\mathbb{N}$ with $|m-n|\le\eta R_n $. If $|m-n|\le\eta R_n$, then \eqref{regularly-varying-R} gives
\[
\frac{m}{n}\in[1-c_0\eta,1+c_0\eta]
\]
for all sufficiently large $n$, where $c_0>0$ is fixed.  The uniform convergence theorem for slowly varying functions \cite{BGT87} gives \(L(m)\asymp L(n)\) uniformly for \(\frac{m}{n}\in[1-c_0\eta,1+c_0\eta]\).  Also $m^{-p}\asymp n^{-p}$. Hence $B_m\asymp B_n$, uniformly for all $m\in\mathbb{N}$ with $|m-n|\le\eta R_n$.
\end{proof}

\begin{prop}\label{stretched-exponential-application}
Fix $0<\alpha<1$ and let
\[
B_n=C_\alpha e^{-n^\alpha},
\qquad
C_\alpha=\Bigl(\sum_{j=1}^{\infty}e^{-j^\alpha}\Bigr)^{-1}.
\]
Then the partition $\beta$ is tail-regular and
\begin{equation}\label{stretched-R}
R_n\asymp n^{1-\alpha}.
\end{equation}
\end{prop}

\begin{proof}
Given $0<\alpha<1$, it is clear that 
\[
\int_n^{\infty}e^{-x^\alpha}\,dx
\le
\sum_{j=n}^{\infty}e^{-j^\alpha}
\le
e^{-n^\alpha}+\int_n^{\infty}e^{-x^\alpha}\,dx, 
\]
and by substituting $u=x^\alpha$, we have 
\[
\int_n^{\infty}e^{-x^\alpha}\,dx
=\frac1\alpha\int_{n^\alpha}^{\infty}u^{\frac1\alpha-1}e^{-u}\,du
\asymp
n^{1-\alpha}e^{-n^\alpha}.
\]
It immediately follows that $B_n^\infty\asymp n^{1-\alpha}B_n.$ Hence   \eqref{stretched-R} holds,  and   $\lim_{n\to\infty}R_n=\infty$.

Next, we prove $\beta$ is tail-regular. Choose $\eta>0$. 
Since $R_n=o(n)$, we choose $n_0\in\mathbb N$ 
such that $\eta R_n\le n/2$ for every $n\ge n_0$. Fix $n\ge n_0$ and $m\in\mathbb N$ with $|m-n|\le\eta R_n$. 
By \eqref{stretched-R}, we have 
\[
|m-n|\lesssim n^{1-\alpha}.
\]
Since $n^{1-\alpha}=o(n)$, we also have $m\asymp n$ for all sufficiently large $n$. By the mean value theorem,
\[
|m^\alpha-n^\alpha|
\le
\alpha\min\{m,n\}^{\alpha-1}|m-n|,
\]
and the right-hand side is uniformly bounded for sufficiently large $n$. Hence
\[
\frac{B_m}{B_n}=e^{-(m^\alpha-n^\alpha)}\asymp1
\]
uniformly whenever $|m-n|\le\eta R_n$. Combining with $\lim_{n\to\infty}R_n=\infty$, we obtain that  $\beta$ is tail-regular.
\end{proof}

\end{document}